\UseRawInputEncoding
\documentclass[12pt,reqno]{amsart}
\usepackage{amsmath,amsfonts,amsthm,amsopn,amssymb}
\usepackage{cite,marginnote}
\usepackage{color,enumitem,graphicx}
\usepackage[colorlinks=true,urlcolor=blue,
citecolor=red,linkcolor=blue,linktocpage,pdfpagelabels,
bookmarksnumbered,bookmarksopen]{hyperref}
\usepackage[english]{babel}

\usepackage[left=2.9cm,right=2.9cm,top=2.8cm,bottom=2.8cm]{geometry}
\usepackage[hyperpageref]{backref}

\numberwithin{equation}{section}

\makeindex
\makeindex
\newtheorem{theorem}{Theorem}[section]
\newtheorem{definition}[theorem]{Definition}
\newtheorem{lemma}[theorem]{Lemma}
\newtheorem{corollary}[theorem]{Corollary}
\newtheorem{proposition}[theorem]{Proposition}
\newtheorem{remark}[theorem]{Remark}

\newcommand{\s}{\section}

\newcommand{\R}{\mathbb R}

\newcommand{\C}{\mathbb C}

\newcommand{\lab}{\label}
\newcommand{\bt}{\begin{theorem}}
\newcommand{\et}{\end{theorem}}
\newcommand{\bl}{\begin{lemma}}
\newcommand{\el}{\end{lemma}}
\newcommand{\bd}{\begin{definition}}
\newcommand{\ed}{\end{definition}}
\newcommand{\bc}{\begin{corollary}}
\newcommand{\ec}{\end{corollary}}
\newcommand{\bp}{\begin{proof}}
\newcommand{\ep}{\end{proof}}
\newcommand{\bx}{\begin{example}}
\newcommand{\ex}{\end{example}}
\newcommand{\bi}{\begin{exercise}}
\newcommand{\ei}{\end{exercise}}
\newcommand{\bo}{\begin{proposition}}
\newcommand{\eo}{\end{proposition}}
\newcommand{\br}{\begin{remark}}
\newcommand{\er}{\end{remark}}
\newcommand{\beq}{\begin{equation}}
\newcommand{\eeq}{\end{equation}}
\newcommand{\ba}{\begin{align}}
\newcommand{\ea}{\end{align}}
\newcommand{\bn}{\begin{enumerate}}
\newcommand{\en}{\end{enumerate}}
\newcommand{\bg}{\begin{align*}}
\newcommand{\bcs}{\begin{cases}}
\newcommand{\ecs}{\end{cases}}

\newcommand{\bean}{\begin{eqnarray*}}
\newcommand{\eean}{\end{eqnarray*}}

\def\C{\mathbb{C}}
\def\N{\mathbb{N}}

\def\R{\mathbb{R}}

\def\bd{\mathrm{bd}\,}

\title[Normalized solution]{Normalized solution to the Sch\"odinger equation with potential and general nonlinear term: Mass super-critical case}

\author[Y.~H.~Ding]{Yanheng Ding}
\author[X.~X.~Zhong]{Xuexiu Zhong}

\address[Y.~H.~Ding]{\newline\indent Institute of Mathematics,Academy of Mathematics and Systems Science
\newline\indent
Chinese Academy of Sciences
\newline\indent
Beijing 100190, PR China}
\email{\href{mailto:dingyh@math.ac.cn}{dingyh@math.ac.cn}}

\address[X.~X.~Zhong]{\newline\indent South China Research Center for Applied Mathematics and Interdisciplinary Studies
\newline\indent
South China Normal University
\newline\indent
Guangzhou 510631, PR China}
\email{\href{mailto:zhongxuexiu1989@163.com}{zhongxuexiu1989@163.com}}

\thanks{Xuexiu Zhong was supported by the NSFC (No.11801581), Guangdong Basic and Applied Basic Research Foundation (2021A1515010034),Guangzhou Basic and Applied Basic Research Foundation(202102020225),Province Natural Science Fund of Guangdong (2018A030310082). }

\subjclass[2000]{}

\begin{document}

\begin{abstract}
In present paper, we prove the existence of solutions $(\lambda, u)\in \R\times H^1(\R^N)$ to the following  Schr\"odinger equation
$$
\begin{cases}
-\Delta u(x)+V(x)u(x)+\lambda u(x)=g(u(x))\quad &\hbox{in}~\R^N\\
0\leq u(x)\in H^1(\R^N), N\geq 3
\end{cases}
$$
satisfying the normalization constraint $\displaystyle \int_{\R^N}u^2 dx=a$. We treat the so-called mass super-critical case here. Under an explicit smallness assumption on $V$ and some Ambrosetti-Rabinowitz type conditions on $g$, we can prove the existence of ground state normalized solutions for prescribed mass $a>0$.
Furthermore, we emphasize that  the mountain pass characterization of a minimizing solution of the problem
$$\inf\left\{\int \left[\frac{1}{2}|\nabla u|^2+\frac{1}{2}V(x)u^2-G(u)\right]dx : \|u\|_{L^2(\R^N)}^{2}=a, P[u]=0\right\},$$
where $G(s)=\int_0^s g(\tau)d\tau$ and
$$P[u]=\int\left[|\nabla u|^2-\frac{1}{2}\langle \nabla V(x), x\rangle u^2 -N\left(\frac{1}{2}g(u)u-G(u)\right)\right]dx.$$

\end{abstract}

\maketitle

\s{Introduction}
\renewcommand{\theequation}{1.\arabic{equation}}
This paper is concerned with the existence of nonnegative solution for some semi-linear elliptic equations in $\R^N$ with prescribed mass ($L^2$-norm). Such problems are motivated in particular by searching for solitary waves (stationary states) in nonlinear equations of the Schr\"odinger or Klei-Gordon type. For concreteness, consider the following nonlinear Schr\"odinger equation
\beq\lab{eq:NLSE}
\begin{cases}
&-i\frac{\partial}{\partial t}\Phi=\Delta \Phi -V(x)\Phi+f(|\Phi|^2)\Phi=0, (x,t)\in \R^N\times \R,\\
&\Phi=\Phi(x,t)\in \C, N\geq 1.
\end{cases}
\eeq
Then looking for a solitary wave in \eqref{eq:NLSE} of ``stationary wave" type, i.e., $\Phi(x,t)=e^{i\lambda t}u(x),\lambda\in \R$, and $u:\R^N\rightarrow \R$, one is led to the equation
\beq\lab{eq:NLSES}
-\Delta u+\left(V(x)+\lambda\right) u=g(u)\;\hbox{in}\;\R^N,
\eeq
where $g(u)=f(|u|^2)u$.
Equation \eqref{eq:NLSES} corresponds to the energy (Lagrangian) functional
$$J_\lambda[u]:=\frac{1}{2}\int_{\R^N}|\nabla u|^2dx+\frac{1}{2}\int_{\R^N}\left(V(x)+\lambda\right)u^2dx-\int_{\R^N}G(u)dx$$
where $G(s):=\int_0^s g(\tau)d\tau$ for $s\in \R$.
For physical reasons, one wants the Lagrangian to be finite, and hence one requires $u$
 to vanish at infinity. This plays the role of a boundary condition for \eqref{eq:NLSES}.
 Therefore, it is usually required that $u\in H^1(\R^N)$.

If fix $\lambda\in \R$ and search for solutions $u$ of \eqref{eq:NLSES}, we call \eqref{eq:NLSES} the {\it fixed frequency problem}. One can apply the variational method, looking for critical points of $J_\lambda[u]$, or some other topological methods, such as fixed point theory, bifurcation or the Lyapunov-Schmidt reduction. The {\it fixed frequency problem} has been widely studied for the decades. A large number of literatures are devoted to study the properties of solutions including the existence, non-existence, multiplicity, asymptotic behavior of the solutions. Some are focused on ground sates or positive solutions, and also some are interested in sign-changing solutions, .etc. It seems that it is impossible to summarize it here since the related literatures are huge.

However, an important, and of course well known, feature of \eqref{eq:NLSE} is conservation of mass:
$$\|\Phi(\cdot,t)\|_{L^2(\R^N)}=\|\Phi(\cdot,0)\|_{L^2(\R^N)}, J[\Phi(\cdot,t)]=J[\Phi(\cdot,0)]\;\hbox{for any}\;t\in \R,$$
where $J$ is the energy functional associated with \eqref{eq:NLSE} defined by
$$J[u]=\frac{1}{2}\int_{\R^N}|\nabla u|^2dx+\frac{1}{2}\int_{\R^N}V(x)u^2dx-\int_{\R^N}G(u)dx,$$
for any
$$u\in \mathcal{H}:=\left\{u\in H^1(\R^N) : \left|\int_{\R^N}V(x)u^2 dx\right|<\infty\right\}.$$
The norm of $\mathcal{H}$ is defined by
\beq
\|u\|_{\mathcal{H}}:=\left(\int_{\R^N}[|\nabla u|^2+V(x)u^2+u^2]dx\right)^{\frac{1}{2}},
\eeq
which is equivalent to the usual norm $\|u\|_{H^1(\R^N)}$ under our assumption $(V1)$ in present paper.
For physical reasons, in spite of the physical relevance of normalized solutions. A natural approach to finding solution of \eqref{eq:NLSES} satisfying the normalization constraint
\begin{equation}\label{eq:norm}
 u\in S_a:=\left\{u\in \mathcal{H}: \int_{\R^N}u^2dx=a\right\}
\end{equation}
consists in finding critical points $u\in \mathcal{H}$ of the energy
$J[u]$ under the constraint \eqref{eq:norm}. Then the parameter $\lambda$ appears as Lagrange multiplier. A solution of \eqref{eq:NLSES} satisfying the prescribed mass constraint \eqref{eq:norm}, we call it {\it fixed mass problem}, which is also called normalized solution in many literatures.

\subsection{A brief review of non-potential case}
The question of finding normalized solutions is already interesting for scalar equations and provides features and difficulties that are not present in the {\it fixed frequency problem}. For the non-potential case, consider the problem
\beq\lab{eq:scalar-set}
-\Delta u+\lambda u=g(u), \int_{\R^N}u^2 dx=a>0, u\in H^1(\R^N).
\eeq
When $g(u)=|u|^{p-1}u$ with $p+1\in (2, {2^*})$, one can show that \eqref{eq:scalar-set} possesses a nontrivial solution $u\in H^1(\R^N)$ only if $\lambda>0$. On the other hand, the positive normalized solution of \eqref{eq:scalar-set} can be completely solved by scaling. Indeed, let $U_p$ be the unique positive radial solution to
\beq\lab{eq:def-Up}
  -\Delta u+u=u^{p}\;\hbox{in}\;\R^N;\quad u(x)\rightarrow 0\ \text{ as $|x|\to\infty$;}
\eeq
cf.\ \cite{Kwong1989}.
Setting
$$U_{\lambda,p}(x):=\lambda^{\frac{1}{p-1}}U_p(\sqrt{\lambda}x),$$
one can check that, up to a translation, $U_{\lambda,p}$ is the unique positive radial solution to
$$ -\Delta u+\lambda u=u^{p}\;\hbox{in}\;\R^N;\quad u(x)\rightarrow 0\ \text{ as $|x|\to\infty$;}.$$
A direct computation shows that
\beq\lab{eq:U-scaling}
\|U_{\lambda}\|_{2}^{2}=\lambda^{\frac{4-(p-1)N}{2(p-1)}}\|U\|_{2}^{2}.
\eeq
So we can see that if $p\neq \bar{p}:=1+\frac{4}{N}$, then there exists a  unique $\displaystyle\lambda_a>0$ such that $\displaystyle \|U_{\lambda_a}(x)\|_{2}^{2}=a$.
That is, there exists a  positive normalized solution to \eqref{eq:scalar-set} for any $a>0$ whenever $p\neq 1+\frac{4}{N}$ (and it is unique up to a translation). While for the so-called mass critical case $p+1=2+\frac{4}{N}$, \eqref{eq:scalar-set} has positive normalized solution if and only if $a=\|U_{1+\frac{4}{N}}\|_{2}^{2}$ (with infinitely many solutions    and  $\lambda>0$).

If $g(u)=\sum_{i=1}^{m}a_i|u|^{\sigma_i}u$ with $m>1, a_i>0$ and
$$\begin{cases}0<\sigma_i<\frac{4}{N-2}\quad &\hbox{if}\;N\geq 3,\\
\sigma_i>0\;&\hbox{if}\;N=1,2.  \end{cases}\;\forall i=1,2,\cdots,m.$$
then $g(u)$ is not of homogeneous, and the scaling method dose not work. And thus the existence of normalized solutions becomes nontrivial and many techniques developed for the {\it fixed frequency problem} can not be applied directly.

A series of theories and tools related to {\it fixed frequency problem} have been developed, such as fixed point theory, bifurcation, the Lyapunov-Schmidt reduction, Nehari manifold method, mountain pass theory and many other linking theories. However, for the {\it fixed mass problem}, the normalization constraint \eqref{eq:norm} certainly brings too much trouble in mathematical treatment.  Comparing to the {\it fixed frequency problem}, the {\it fixed mass problem} possesses the following technical difficulties when dealing with it in the variational framework:
\begin{itemize}
\item[(1)] One can not use the usual Nehari manifold method since the frequency is unknown;
\item[(2)] The existence of bounded Palais-Smale sequences requires new arguments;
\item[(3)]  The Lagrange multipliers have to be controlled;
\item[(3)] The embedding $H^1(\R^N)\hookrightarrow L^2(\R^N)$ is not compact. For the {\it fixed frequency problem}, usually a nontrivial weak limit is also a solution. However, for the {\it fixed mass problem}, even the weak limit is nontrivial, the constraint condition may be not satisfied.
\item[(4)] There is a number $\bar{p}=2+\frac{4}{N}$, called the mass critical exponent, affects the geometry of the functional heavily.
\end{itemize}

The best studied cases of \eqref{eq:scalar-set}  correspond to the situation in which a solution can be found as a global minimizer of the energy functional on $S_a$, which is the case of  mass sub-critical for $g(u)$. This research mainly started with the work of Stuart\cite{Stuart1980,Stuart1981,Stuart1989}.  Stuart applied a bifurcation approach to study the nonhomogeneous nonlinearities (not necessarily autonomous). However, it requires more restrictive growth conditions on $g$ to guarantee the compactness. Stated in the particular case of \eqref{eq:scalar-set} with $g(u)=\sum_{i=1}^{m}a_i|u|^{\sigma_i}u$,
 Stuart studied the mass sub-critical case, i.e., $0<\sigma_i<\frac{4}{N},1\leq i\leq m$, he showed that the corresponding functional is bounded from  below  on $S_a$ for any $a>0$ and obtained a critical point proving that the infimum is reached (see \cite{Stuart1981,Stuart1989} for the details).
Later on, the concentration compactness principle of Lions\cite{Lions1984a,Lions1984b} was used to study this kind problem. We own to M.Shibata \cite{Shibata2014} the study of the mass sub-critical case for general nonlinearities.

However, when the nonlinearity $g(u)$ is of mass super-critical, the energy functional  is unbounded below on $S_a$. Thus it is impossible to search for a minimum on $S_a$.  L. Jeanjean \cite{Jeanjean1997} could prove the corresponding functional possesses the mountain pass geometric structure on $S_a$. Then obtained  the normalized solution by a minimax approach and a smart compactness argument.

 A multiplicity result was established by Bartsch and de Valeriola in \cite{bartsch2003multi}.  Ikoma and Tanaka  provided an alternative proof for this multiplicity result by exploiting an idea related to symmetric mountain pass theorems, see \cite{ikomatanaka2019}. The more general nonlinearities are studied in \cite{sslu2020cvpde,Mederski2020}.

\subsection{A brief review of the potential case}
Observing that in practice physical background, many problems inevitably involve the potential $V(x)\not\equiv 0$. Recently, the mass prescribed  problem with potential
\beq\lab{eq:V-pontential}
-\Delta u+(V(x)+\lambda)u=g(u)\;\hbox{in}~\R^N, \int_{\R^N}u^2 dx=a>0
\eeq
is also studied under some suitable assumptions on $V(x)$ by many researchers.
We note that the techniques developed in a sequences of literature mentioned above study the non-potential case, can not applied directly. Hence, besides the importance in the applications, not negligible reasons of many researchers' interest for such problems are its stimulating and challenging difficulties.

In \cite[Section 3]{PellacciPistoiaVairaVerzini}, Pellacci et al. apply Lyapunov-Schmidt reduction approach to study the special case of $g(u)=u^p$ and some existence results are obtained. When dealing with the general nonlinearities of mass-subcritical case, the main difficulty is the compactness of the minimizing sequence. And the key step is to establish the strict so-called sub-additive inequality. In \cite{Ikoma-Miyamoto2020}, Norihisa Ikoma and Yasuhito Miyamoto applied the standard concentration compactness arguments, which is due to Lions\cite{Lions1984a,Lions1984b}, to make the first progress on this direction. Recently, under different assumptions, Zhong and Zou \cite{zhongzou2020} present a new approach based on iteration, to establish the strict sub-additive inequality, which can greatly simplify the discussion process in the traditional sense.

While dealing with the mass super-critical case involves  potential, the functional is unbounded from below, thus one can not apply the minimizing argument constraint on the sphere $S_a$ any more. The problem becomes more complicated.
The case of positive potential $V(x)\geq 0$ and vanishing at infinity, including potentials with singularities, is considered in the very recent paper \cite{Bartsch2021}. In such a case, the mountain pass structure by Jeanjean\cite{Jeanjean1997} is destroyed, but a new variational principle exploiting the Pohozaev identity can be used.  By constructing a suitable linking structure, the authors of \cite{Bartsch2021} can obtain the existence of solutions with high Morse index.
The case of negative potential $V(x)\leq 0$ and vanishing at infinity, is also considered in \cite{MG2021}. Under some explicit smallness assumption on $V(x)$, the authors of \cite{MG2021} can obtain the existence of solutions with prescribed $L^2$-norm.
However, for the the mass super-critical case involves  potential, all the literatures mentioned above,  only consider the case of $g(u)=|u|^{p-2}u$. And it seems that their discussion can not be directly extended to general nonlinear terms.

In present paper, we shall treat the mass super-critical case with negative potential and general nonlinear terms. Precisely, we consider
\beq
g(u)=|u|^{p-2}u+h(u),
\eeq
satisfying some Ambrosetti-Rabinowitz type condition. We point out that in our argument, $g(u)$ is not a  perturbation of $|u|^{p-2}u$. Indeed, under some suitable assumptions on $V(x)$, it is easy to check that the mountain pass structure by Jeanjean\cite{Jeanjean1997} is kept. However, since $V(x)$ is not necessary radial, we can not work in the radial sub-space $H_r^1(\R^N)$. Hence, the main difficulty is the compactness. Under some suitable assumptions on $V(x)$ and $g(u)$, here we can prove the existence of ground state normalized solution as well as mountain pass characterization of a minimizing solution of the problem
\beq\lab{eq:mini-problem}
\inf\left\{J[u]: u\in S_a, P[u]=0\right\}, G(\xi)=\int_0^\xi g(\tau)d\tau,
\eeq
and
\beq\lab{eq:P-identity}
P[u]:=\|\nabla u\|_2^2-\frac{1}{2}\int_{\R^N}\langle\nabla V(x), x\rangle u^2dx-N\int_{\R^N}\left[\frac{1}{2}g(u)u-G(u)\right]dx.
\eeq

\br\lab{remark:r1}
\begin{itemize}
\item[(1)]When study the mass super-critical problem involving potential, we emphasize that this is the first paper  considering the general nonlinearities.
\item[(2)]
When the nonlinearities are combination of  mass-subcritical and mass-supercritical terms, the geometric structure of the corresponding functional will become very complex.  And  ones  need to give more detailed  arguments on the geometry of the corresponding functional constrained on suitable sub-manifold. For the non-potential case, we refer to \cite{soave2020,soave2020cri}. And we note that so far there is no literature about the  combination case with potential.
\item[(3)]We also remark that there are also some works on the system problem with potentials. Consider the following system:
    \beq\lab{eq:20210414-system}
\begin{cases}
-\Delta u_1+V_1(x)u_1+\lambda_1 u_1=\mu_1|u_1|^{p_1-2}u_1+ \beta r_1|u_1|^{r_1-2}|u_2|^{r_2} u_1\;\quad&\hbox{in}\;\R^N,\\
-\Delta u_2+V_2(x)u_2+\lambda_2 u_2=\nu_1|u_2|^{q_1-2}u_2+\beta r_2 |u_1|^{r_1}|u_2|^{r_2-2} u_2 \;\quad&\hbox{in}\;\R^N,\\
u_1,u_2\in H^1(\R^N), N\geq 1.
\end{cases}
\eeq
 Ikoma and Miyamoto \cite{Ikoma2021} give a first progress for $r_1=r_2$. And the case $r_1\neq r_2$ is owed to Deng,He and Zhong \cite{DHZ2021}. However all of them only study the mass sub-critical case, which enable them to apply a minimizing argument on the $L^2$-sphere $S_a$. We also note that the case of mass super-critical system with potentials is also blank, even for $g(u)=u^{p-1}$.
\end{itemize}
\er

The paper is organized as follows. In the next section we shall propose the hypotheses and discuss our results. In section \ref{sec:preliminaries}, we study the related properties of the so-called Pohozaev manifold, and the mini-max structure. In Section \ref{Proof-th1}, we prove our main Theorem.
Throughout the paper we use the notation $\|u\|_p$ to denote the $L^p$-norm. The notation $\rightharpoonup$ denotes weak convergence in $H^1(\R^N)$. Capital latter $C$ stands for positive constant which may depend on some parameters, whose precise value can change from line to line.
\medskip

\s{Hypotheses and statement of results}\lab{sec:Hypotheses-statements}
\renewcommand{\theequation}{2.\arabic{equation}}
For the nonlinearities, we assume the following hold.
\begin{itemize}
\item[(G1)] $g:\R\rightarrow \R$ is continuous and odd.
\item[(G2)]There exists some $(\alpha,\beta)\in \R_+^2$ satisfying $2+\frac{4}{N}<\alpha\leq \beta<\frac{2N}{N-2}$ such that
\beq\lab{eq:condition-G1}
\alpha G(s)\leq g(s)s\leq \beta G(s)~\hbox{with}~G(s)=\int_0^s g(t)dt.
\eeq
\item[(G3)] The functional defined by $\widetilde{G}(s):=\frac{1}{2}g(s)s-G(s)$ is of class $C^1$ and
    $$\widetilde{G}'(s)s\geq \alpha \widetilde{G}(s),\forall s\in \R,$$
where $\alpha$ is given by (G2).
\end{itemize}

For the non-potential case, i.e., $V(x)\equiv 0$, under the assumptions $(G1)$ and $(G2)$,
Jeanjean \cite{Jeanjean1997} can prove that $J[u]$ satisfies mountain pass geometry and thus a bounded Palais-Smale sequence exists.
Recalling that the embedding $H_r^1(\R^N)\hookrightarrow L^p(\R^N)$ is compact for any $p\in (2, 2^*)$ if $N\geq 2$, Jeanjean worked in the radial space $H_r^1(\R^N)$ and obtain the normalized solution for any $a>0$ provided $N\geq 2$ (see \cite[Theorem 2.1]{Jeanjean1997}). And if furthermore the following (G'3) holds,
\begin{itemize}
\item[(G'3)] $\widetilde{G}'(s)s>\left(2+\frac{4}{N}\right)\widetilde{G}(s),\forall s\in \R\backslash\{0\}$.
\end{itemize}
the result will be valid for all $N\geq 1$, see \cite[Theorem 2.2]{Jeanjean1997}.

Recalling the fiber map
\beq\lab{eq:def-fiber-map}
u(x)\mapsto (t\star u)(x):=t^{\frac{N}{2}}u(tx),
\eeq
for $(t, u)\in \R^+\times H^1(\R^N)$, which preserves the $L^2$-norm. Define
\beq\lab{eq:20210529-e1}
 \Psi_{u}(t)=J[t\star u]~\hbox{and}~\Psi_{\infty,u}(t)=I[t\star u]
\eeq
with
\beq\lab{eq:def-I}
I[u]=\frac{1}{2}\int_{\R^N}|\nabla u|^2dx-\int_{\R^N}G(u)dx.
\eeq
A direct computation shows that
\beq\lab{eq:20210529-e2}
(\Psi_{u})'(t)=\frac{1}{t}P[t\star u]~\hbox{and}~(\Psi_{\infty,u})'(t)=\frac{1}{t}P_\infty[t\star u]
\eeq
where
\beq\lab{eq:def-P}
P[u]:=\|\nabla u\|_2^2-\frac{1}{2}\int_{\R^N}\langle \nabla V(x), x\rangle u^2 dx-N\int_{\R^N}\left[\frac{1}{2}g(u)u-G(u)\right]dx.
\eeq
and
\beq\lab{eq:def-Pinfty}
P_\infty[u]:=\|\nabla u\|_2^2-N\int_{\R^N}\left[\frac{1}{2}g(u)u-G(u)\right]dx.
\eeq
Define
\beq\lab{eq:def-manifold-P-and-Pinfty}
\mathcal{P}:=\{u\in H^1(\R^N):P[u]=0\}\;\hbox{and}~\mathcal{P}_\infty:=\{u\in H^1(\R^N):P_\infty[u]=0\}
\eeq
For $a>0$, set
\beq\lab{eq:20210529-e3}
\mathcal{P}_{a}:=S_a\cap \mathcal{P}~\hbox{and}~\mathcal{P}_{\infty,a}:=S_a\cap \mathcal{P}_\infty.
\eeq
\br\lab{remark:for-I}
\begin{itemize}
\item[(1)]
We remark that (G'3) plays an crucial role to guarantee an unique $t_u>0$ such that $t_u\star u\in \mathcal{P}_{\infty, a}$ for any $u\in S_a$(see \cite[Lemma 2.9]{Jeanjean1997}).  Furthermore,
\beq\lab{eq:20210609-ze1}
I[t_u\star u]=\max_{t>0} I[t\star u].
\eeq
And thus
\beq\lab{eq:def-ma-20210615}
m_a:=\inf_{u\in S_a}\max_{t>0}I[t\star u]=\inf_{u\in \mathcal{P}_{\infty,a}}I[u].
\eeq
So one can see that the mountain pass given by Jeanjean \cite[Theorem 2.2]{Jeanjean1997} characterization of a minimizing solution of the problem
$\displaystyle\inf_{u\in \mathcal{P}_{\infty,a}}I[u]$.
\item[(2)] Under the assumptions (G1),(G2) and (G'3), the uniqueness $t_u$ can imply the strict decreasing of $m_a$ for $N\geq 2$, see\cite[Theorem 1.1]{Yang2020}. Then Yang applied a minimizing argument to give a new proof to \cite[Theorem 2.2]{Jeanjean1997}.
\item[(3)]We remark that to obtain the convergence of the Palais-Smale sequence, the author of \cite{Yang2020} worked in the subspace $H_r^1(\R^N)$ since that the equation is autonomous.
\end{itemize}
\er

In present paper, we study the potential case $0\not\equiv V(x)$. We emphasize that $V(x)$ is not necessary radial. Hence, we can not work in the radial subspace $H_r^1(\R^N)$. For the trapping potentials, i.e.,
$\displaystyle\lim_{|x|\rightarrow +\infty}V(x)=+\infty,$
we remark that this coercive condition can provide enough compactness to cause the existence of solutions which are local minimizers of $J$ on $S_a$ also in the mass-supercritical case, at least when $a$ is sufficiently small, see \cite{BBJV2017,NTV2015}. Hence, in present paper we only consider the non-trapping potentials case. Suppose that $\displaystyle\limsup_{|x|\rightarrow +\infty}V(x)=:V_\infty<+\infty$ exists. And without loss of generality, we may assume that $V_{\infty}=0$. If not, since we are finding the normalized solution, we may replace $(V(x), \lambda)$ by
$\displaystyle (\tilde{V}(x), \tilde{\lambda}):=(V(x)-V_{\infty}, \lambda+V_{\infty})$.

Precisely, we assume that
\begin{itemize}
\item[$(V1)$] $\displaystyle\lim_{|x|\rightarrow +\infty} V(x)=\sup_{x\in \R^N}V(x)=0$ and there exists some $\sigma_1\in \left[0,\frac{N(\alpha-2)-4}{N(\alpha-2)}\right)$ such that
\beq
\left|\int_{\R^N}V(x)u^2 dx\right|\leq \sigma_1 \|\nabla u\|_2^2, \forall u\in H^1(\R^N).
\eeq
\item[$(V2)$] $\nabla V(x)$ exists for a.e. $x\in \R^N$, put $W(x):=\frac{1}{2}\langle \nabla V(x),x\rangle$. There exists some $0\leq \sigma_2<\min\{\frac{N(\alpha-2)(1-\sigma_1)}{4}-1, \frac{N}{\beta}-\frac{N-2}{2}\}$ such that
\beq
\left|\int_{\R^N}W(x)u^2 dx\right|\leq \sigma_2 \|\nabla u\|_2^2, \forall u\in H^1(\R^N).
\eeq
\item[$(V3)$] $\nabla W(x)$ exists for a.e. $x\in \R^N$, put
$$Y(x):=(\frac{N}{2}\alpha-N)W(x)+\langle \nabla W(x),x\rangle.$$
There exists some $\sigma_3\in [0,\frac{N}{2}\alpha-N-2)$ such that
\beq
\int_{\R^N} Y_+(x)u^2 dx\leq \sigma_3\|\nabla u\|_2^2,~\forall u\in H^1(\R^N).
\eeq
\end{itemize}

\br
For $N\geq 3$, by Sobolev inequality, under some small conditions on $\|V\|_{\frac{N}{2}},\|W\|_{\frac{N}{2}}$ and $\|Y_+\|_{\frac{N}{2}}$, then $(V1),(V2)$ and $(V3)$ hold. For example, let $S$ be the sharp constant in the critical Sobolev inequality, if $S\|V\|_{\frac{N}{2}}<\frac{N(\alpha-2)-4}{N(\alpha-2)}$, then we can take $\sigma_1:=\|V\|_{\frac{N}{2}}S$ and
\begin{align*}
\left|\int_{\R^N}V(x)u^2 dx\right|\leq \|V\|_{\frac{N}{2}} \|u\|_{2^*}^{2}\leq S\|V\|_{\frac{N}{2}} \|\nabla u\|_2^2=\sigma_1 \|\nabla u\|_2^2.
\end{align*}
Another important application is the Hardy potentials $V(x)=-\frac{\mu}{|x|^2}$ with suitable small $\mu$.
\er

We are aim to establish the following result.
\bt\lab{thm:non-potential}
Let $N\geq 3$.
Assume that $g(s)$ satisfies the hypotheses (G1)-(G3) and $0\not\equiv V(x)$ satisfies the conditions $(V1)$-$(V3)$. Then for any $a>0$,
there exists a couple $(\lambda_a, u_a)\in \R^+\times H^1(\R^N)$ solves
\beq\lab{eq:problem-nonpotential}
-\Delta u(x)+V(x)u(x)+\lambda u(x)=g(u(x)), \;\lambda\in \R, x\in\R^N,
\eeq
satisfying $\displaystyle \int_{\R^N}u^2 dx=a$.
\et

\br
In a sequence of literatures, one often need $N\leq 4$ to apply the Liouvill type result \cite[Lemma A.2]{ikoma2014compactness} which indicates that  if a function $u(x)$ satisfying
$$-\Delta u\geq 0\;\hbox{in}\;\R^N \;\hbox{and}\;u\in L^p(\R^N) \;\hbox{with}~\begin{cases}p<+\infty, N=1,2,\\
p\leq \frac{N}{N-2}, N\geq 3 \end{cases}$$
then $u(x)\equiv 0$ in $\R^N$. Since $u\in L^2(\R^N)$, for the case of $N\leq 4$, one can prove that $\lambda>0$ prior to the compactness. In present paper,  our Theorem \ref{thm:non-potential} is also valid for $N\geq 5$.
\er

\br
For the Lagrange multiplier $\lambda_a$, from a physical point of view it represents the chemical potentials of standing waves. We point out that there are situations in the Bose-Einstein condensate theory that requires the chemical potentials are positive, see \cite{LiebSeiringerSolovejYngvason2005,Pitaevskii2003}. By a well known result own to Frank H. Clarke(see \cite[Theorem 1]{Clarke1976}), one can also always have that $\lambda_a\geq 0$. From a purely mathematical point of view, under the Ambrosetti-Rabinowitz type condition, this conclusion is more intuitive, see Lemma \ref{lemma:non-existence} below.
\er

\br\lab{remark:20210615-r1}
In the setting of Theorem \ref{thm:non-potential}, a solution is obtained through a mini-max method:
\beq\lab{eq:def-Ca-20210615}
C_a:=\inf_{u\in S_a}\max_{t>0} J[t\star u].
\eeq
 One can also easy to see that $C_a=\inf_{u\in \mathcal{P}_a}J[u]=\gamma_a$, where $\gamma_a$ is the mountain pass value defined by Jeanjean\cite[Lemma 2.8]{Jeanjean1997}. See \cite[Lemma 2.10]{Jeanjean1997}.
\er

\s{Pohozaev constraint and a mini-max structure}\lab{sec:preliminaries}
\renewcommand{\theequation}{3.\arabic{equation}}
\bl\lab{lemma:Pohozaev}
Assume that $u\in H^1(\R^N)$ is a solution to \eqref{eq:NLSES},then
$u\in \mathcal{P}$.
\el
\bp
Let $u$ be a solution to \eqref{eq:NLSES}, then testing \eqref{eq:NLSES} by $u$, we have
\beq\lab{eq:Nehari}
\int_{\R^N}|\nabla u|^2dx+\int_{\R^N}(V(x)+\lambda)u^2dx -\int_{\R^N}g(u)udx=0.
\eeq
On the other hand, a solution to \eqref{eq:NLSES} must satisfy the so-called Pohozaev identity
\begin{align}\lab{eq:Pohozaev}
&(N-2)\int_{\R^N}|\nabla u|^2dx+N\int_{\R^N}(V(x)+\lambda)u^2dx  \nonumber \\
&-2N\int_{\R^N}G(u)dx+\int_{\R^N}\langle \nabla V(x), x\rangle u^2 dx=0.
\end{align}
Eliminate the unknown parameter $\lambda$, we obtain that
$$\int_{\R^N}|\nabla u|^2dx-\int_{\R^N}W(x) u^2 dx-N\int_{\R^N}\left[\frac{1}{2}g(u)u-G(u)\right]dx=0.$$
\ep

\bo\lab{prop:20210529-p1}
Let $u\in S_a$. Then: $t\in \R^+$ is a critical point for $\Psi_{u}(t):=J[t\star u]$ if and only if $t\star u\in \mathcal{P}_{a}$.
\eo
\bp
By a direct computation, it follows by \eqref{eq:20210529-e2} that $\Psi'_u(t)=\frac{1}{t}P[t\star u]$.
\ep

\bo\lab{prop:20210615-p1}
For any critical point of $J\big|_{\mathcal{P}_{a}}$ , if $(\Psi_{u})''(1)\neq 0$, then there exists some $\lambda\in \R$ such that
$$J'[u]+\lambda u=0.$$
\eo
\bp
Let $u$ be a critical point of $J[u]$ constraint on $\mathcal{P}_{a}$, then there exists $\lambda,\mu\in \R$ such that
\beq\lab{eq:20210507-e1}
J'[u]+\lambda u+\mu P'(u)=0.
\eeq
We only need to prove that $\mu=0$.
Noting that a function $u$ solves \eqref{eq:20210507-e1} must satisfy the corresponding Pohozave identity
$$\frac{d}{dt} \Phi[t\star u]\Big|_{t=1}=0,$$
where $\Phi[u]:=J[u]+\frac{1}{2}\lambda \|u\|_2^2+\mu P[u]$ is the corresponding energy functional of \eqref{eq:20210507-e1}.
Observing that
$$\Phi[t\star u]=J[t\star u]+\frac{1}{2}\lambda \|u\|_2^2+\mu P[t\star u]
=\Psi_{u}(t)+\frac{1}{2}\lambda \|u\|_2^2+\mu \frac{1}{t}(\Psi_{u})'(t),$$
we have
\begin{align*}
\frac{d}{dt} \Phi[t\star u]=& (1-\frac{\mu}{t^2})(\Psi_{u})'(t)+\mu \frac{1}{t}(\Psi_{u})''(t).
\end{align*}
Hence,
\begin{align*}
0=&\frac{d}{dt} \Phi[t\star u]\Big|_{t=1}=(1-\mu)(\Psi_{u})'(1)+\mu (\Psi_{u})''(1)\\
=&(1-\mu)P[u]+\mu ~(\Psi_{u})''(1)=\mu ~(\Psi_{u})''(1).
\end{align*}
So by $(\Psi_{u})''(1)\neq 0$, we have that $\mu=0$.
\ep

\bl\lab{lemma:non-existence}
 Assume that there exists some $2<\alpha< \beta<2^*$ such that
\beq\lab{eq:AR-condition}
\alpha G(s)\leq g(s)s\leq \beta G(s),~\forall s\in \R.
\eeq
Then the following equation
\beq\lab{eq:limit-equation-1}
-\Delta u+\lambda u=g(u)
\eeq
has no nontrivial  solution $u\in H^1(\R^N)$ provided $\lambda\leq 0$.
\el
\bp
We argue by contradiction and suppose that there exists some  $0\not\equiv u\in H^1(\R^N)$ solves \eqref{eq:limit-equation-1} with $\lambda\leq 0$. Testing \eqref{eq:limit-equation-1} by $u$, we get that
\beq\lab{eq:20210506-e1}
\|\nabla u\|_2^2+\lambda \|u\|_2^2=\int_{\R^N}g(u)u dx.
\eeq
On the other hand, since $u$ is a solution to  \eqref{eq:limit-equation-1}, we have that
\beq\lab{eq:20210506-e2}
(N-2)\|\nabla u\|_2^2+N\lambda\|u\|_2^2 -2N\int_{\R^N}G(u)dx=0.
\eeq
Then by \eqref{eq:20210506-e1} and \eqref{eq:20210506-e2}, we can find some $C>0$ such that
\beq\lab{eq:20210506-e3}
\lambda \|u\|_2^2 =\int_{\R^N} \left[N G(u)-\frac{N-2}{2}g(u)u\right] dx\geq C\int_{\R^N}G(u)dx .
\eeq
Indeed, under the assumptions, if $N=1,2$, it is easy to see that
$$\int_{\R^N} \left[N G(u)-\frac{N-2}{2}g(u)u\right] dx\geq N\int_{\R^N}G(u)dx.$$
 And for the case of $N\geq 3$, we also have that
 $$\int_{\R^N} \left[N G(u)-\frac{N-2}{2}g(u)u\right] dx\geq \int_{\R^N}[NG(u)-\frac{N-2}{2}\beta G(u)]dx=C\int_{\R^N}G(u)dx$$
 with $C:=N-\frac{N-2}{2}\beta>0$ due to the fact $\beta<\frac{2N}{N-2}$.
 Hence, if $\lambda<0$, by $0\not\equiv u\in H^1(\R^N)$ and \eqref{eq:20210506-e3}, we have that
 $$0>\lambda \|u\|_2^2\geq C \int_{\R^N}G(u)dx\geq 0,$$
 a contradiction.
 If $\lambda=0$, then by \eqref{eq:20210506-e3} again, have $\int_{\R^N}G(u)dx=0$. Then it follows \eqref{eq:AR-condition} that $\int_{\R^N}g(u)u=0$. Hence, by \eqref{eq:20210506-e1}, we obtain that $\|\nabla u\|_2^2=0$, a contradiction to $0\not\equiv u\in H^1(\R^N)$.
\ep

\bl\lab{lemma:20210611-zl1}
Assume that the hypotheses (G1) and (G2) hold.
Then for any $a>0$, there exists some $\delta_a>0$ such that
\beq\lab{eq:20210611-e1}
\inf\left\{t>0:\exists u\in S_a~\hbox{with}~\|\nabla u\|_2^2=1~\hbox{such that}~t\star u\in \mathcal{P}_{a}\right\}\geq \delta_a.
\eeq
Consequently,
\beq\lab{eq:20210926-xe1}
\inf_{u\in \mathcal{P}_a} \|\nabla u\|_2^2 \geq \delta_a^2.
\eeq
\el
\bp
 A direct computation shows that
\begin{align*}
(\Psi_{u})'(t)=&\|\nabla u\|_2^2 t-\int_{\R^N}W(\frac{x}{t})u^2dx t^{-1}-
N\int_{\R^N}\widetilde{G}(t^{\frac{N}{2}}u(x))dx t^{-N-1}.
\end{align*}
Under the assumption $(V2)$, for any $u\in D_{0}^{1,2}(\R^N)$, one can see that
\begin{align*}
&\|\nabla u\|_2^2 t-\int_{\R^N}W(\frac{x}{t})u^2dx t^{-1}\\
=&\|\nabla u\|_2^2 t-\int_{\R^N}W(x)(t\star u)^2dx t^{-1}\\
\geq& (1-\sigma_2) \|\nabla u\|_2^2 t.
\end{align*}
So that for any $u\in S_a$ with $\|\nabla u\|_2^2=1$, if $t\star u\in \mathcal{P}_{a}$, by
Proposition \ref{prop:20210529-p1}, we have $(\Psi_{u})'(t)=0$ and thus
\begin{align*}
(1-\sigma_2)\leq&\|\nabla u\|_2^2 -\int_{\R^N}W(\frac{x}{t})u^2dx t^{-2}\\
=&N\int_{\R^N}\widetilde{G}(t^{\frac{N}{2}}u(x))dx t^{-N-2}
\end{align*}
That is,
\beq\lab{eq:20210611-e3}
1-\sigma_2\leq\frac{N}{2}\int_{\R^N}g(t^{\frac{N}{2}}u(x))
u(x)t^{-\frac{N}{2}-2}dx-N\int_{\R^N}G(t^{\frac{N}{2}}u(x))t^{-N-2}dx.
\eeq
By (G2),
\beq\lab{eq:20210611-e4}
\frac{1}{\beta}g(t^{\frac{N}{2}}u(x))t^{\frac{N}{2}}u(x)\leq G(t^{\frac{N}{2}}u(x)).
\eeq
So by \eqref{eq:20210611-e3} and \eqref{eq:20210611-e4}, we obtain that
\beq\lab{eq:20210611-e5}
 1-\sigma_2\leq N\left(\frac{1}{2}-\frac{1}{\beta}\right) t^{-N-2} \int_{\R^N}g(t^{\frac{N}{2}}u(x)) t^{\frac{N}{2}}u(x) dx.
\eeq
By (G2) again, one can find some suitable $C_1>0$ such that
$g(s)s\leq C_1(s^\alpha+s^\beta), \forall s\in \R$.
On the other hand, by Gagliardo-Nirenberg  inequality, we can find some $C_2>0$ independent of the choice of $u$, such that
\beq\lab{eq:20210611-e6-bu}
\|u\|_\alpha^\alpha\leq C_2, \|u\|_\beta^\beta\leq C_2.
\eeq
And thus
\beq\lab{eq:20210611-e6}
1-\sigma_2\leq C_1C_2N\left(\frac{1}{2}-\frac{1}{\beta}\right)\left[t^{\frac{N}{2}\alpha-N-2}+t^{\frac{N}{2}\beta-N-2}\right].
\eeq
So by $2+\frac{4}{N}<\alpha<\beta$, we obtain the lower bound $\delta_a>0$.
\ep
Consider the decomposition of $\mathcal{P}_{a}$ into the disjoint union
$$\mathcal{P}_{a}=\mathcal{P}_{a}^{+}\cup \mathcal{P}_{a}^{0}\cup \mathcal{P}_{a}^{-},$$
where
\beq
\mathcal{P}_{a}^{+(resp.~0,-)}:=\left\{u\in \mathcal{P}_{a}: (\Psi_{u})''(1)>(resp.~=,<)0\right\}.
\eeq

\bl\lab{lemma:20210609-xl1}
Assume that the hypotheses (G1)-(G3) hold, suppose further $(V3)$, then $\mathcal{P}_{a}^{-}=\mathcal{P}_{a}$ is closed in $H^1(\R^N)$ and it is a natural constraint of $J\big|_{S_a}$
\el
\bp
Recalling $W(x)=\frac{1}{2}\langle \nabla V(x), x\rangle$,
for any $u\in \mathcal{P}$, we have that
\beq\lab{eq:20210609-zdle3}
\int_{\R^N}|\nabla u|^2dx-\int_{\R^N}W(x) u^2 dx=N\int_{\R^N}\widetilde{G}(u)dx.
\eeq
A direct computation shows that
\begin{align*}
(\Psi_{u})''(1)=&\|\nabla u\|_2^2+\int_{\R^N}W(x)u^2dx+\int_{\R^N}\langle \nabla W(x),x\rangle u^2dx\\
&+N(N+1)\int_{\R^N}\widetilde{G}(u)dx-\frac{N^2}{2}\int_{\R^N}\widetilde{G}'(u)udx.
\end{align*}
Using (G3), $(V3)$ and \eqref{eq:20210609-zdle3}, we have that
\begin{align}\lab{eq:20210616-e1}
(\Psi_{u})''(1)\leq &\int_{\R^N}Y(x)u^2 dx-\left(\frac{N}{2}\alpha-N-2\right)\|\nabla u\|_2^2\nonumber\\
\leq& -\left(\frac{N}{2}\alpha-N-2-\sigma_3\right)\|\nabla u\|_2^2<0.
\end{align}
So that $\mathcal{P}_{a}^{+}=\mathcal{P}_{a}^{0}=\emptyset$, which implies that
$\mathcal{P}_{a}^{-}=\mathcal{P}_{a}$ is closed in $H^1(\R^N)$. Furthermore, by Proposition \ref{prop:20210615-p1}, one can see that it is a natural constraint of $J\big|_{S_a}$.
\ep

\br
Let $\{u_n\}\subset \mathcal{P}_{a}^{-}$ such that $J[u_n]$ approaches a possible critical value. Noting that $(\Psi_{u_n})''(1)<0$ is an open constraint, then there exists sequences $\lambda_n,\mu_n\in \R$ such that
$$J'[u_n]+\lambda_n u_n +\mu_n \mathcal{P}'[u_n]\rightarrow 0.$$
Applying a similar argument as Proposition \ref{prop:20210615-p1}, we obtain that
\beq\lab{bu-hll-e2}
\mu_n (\Psi_{u_n})''(1)\rightarrow 0.
\eeq
By \eqref{eq:20210616-e1}, one can see that
$$\mu_n \|\nabla u_n\|_2^2\rightarrow 0.$$
However, by Lemma \ref{lemma:20210611-zl1},  we have that $\|\nabla u\|_2^2\geq \delta_a^2>0$ for all $u\in \mathcal{P}_a$. Hence, we obtain that $\mu_n\rightarrow 0$.
And thus if furthermore $\{u_n\}$ is bounded in $H^1(\R^N)$, then we have that
$$J'[u_n]+\lambda_n u_n\rightarrow 0\;\hbox{in}\;H^{-1}(\R^N).$$
\er

\bc\lab{cro:20210616-c1}
Under the assumptions (G1)-(G3) and $(V1)$-$(V3)$, for any $u\in H^1(\R^N)\backslash\{0\}$, there exists an unique $t_u>0$ such that $t_u\star u\in \mathcal{P}$. Furthermore,
\beq
J[t_u\star u]=\max_{t>0}J[t\star u].
\eeq
\ec
\bp
Let $a:=\|u\|_2^2$, since $u\in H^1(\R^N)$, we have that $\|\nabla u\|_2>0$. Then by a similar proof of Lemma \ref{lemma:20210611-zl1}, one can see that
$\Psi'_u(t)>0$ for $t$ small enough. Hence, there exists some $t_0>0$ such that $\Psi_u(t)$ increases in $t\in (0,t_0)$. On the other hand, under the assumptions (G1)-(G3), by $(V1)$, it is not hard to see that
\beq\lab{eq:20210616-xe1}
\lim_{t\rightarrow +\infty} J[t\star u]=-\infty.
\eeq
So there must exists some $t_1>t_0$ such that
$$J[t_1\star u]=\max_{t>0}J[t\star u].$$
Hence, $\Psi'_u(t_1)=0$ and it follows Proposition \ref{prop:20210529-p1} that $t_1\star u\in \mathcal{P}$.
Suppose that there exists another $t_2>0$ such that $t_2\star u\in \mathcal{P}$. Then by Lemma \ref{lemma:20210609-xl1}, we have that both $t_1$ and $t_2$ are strict local maximum of $\Psi_u(t)$. Without loss of generality, we assume that $t_1<t_2$. Then there exists some $t_3\in (t_1, t_2)$ such that
$$\Psi_u(t_3)=\min_{t\in [t_1,t_2]}\Psi_u(t).$$
So we have that $\Psi'_u(t_3)=0$ and $\Psi''_{u}(t_3)\geq 0$,  which implies that
$t_3\star u \in \mathcal{P}_{a}^{+}\cup \mathcal{P}_{a}^{0}$,
a contradiction to Lemma \ref{lemma:20210609-xl1}.
\ep

\br\lab{remark:Ca-positive}
Under the hypotheses (G1)-(G3) and $(V1)$-$(V3)$, for any $u\in \mathcal{P}_{a}$,
one can see that $\Psi_u(t)\rightarrow 0$ as $t\rightarrow 0^+$ and $\Psi_u(t)\rightarrow -\infty$ as $t\rightarrow +\infty$. By Corollary \ref{cro:20210616-c1}, we have that
$$J[u]=\max_{t>0}J[t\star u]>0.$$
Furthermore,
\beq
C_a:=\inf_{u\in \mathcal{P}_a}J[u]=\inf_{u\in S_a}\max_{t>0}J[t\star u].
\eeq
On the other hand, by (G2) and $(V2)$, we have that
\begin{align*}
(1+\sigma_2)\|\nabla u\|_2^2\geq & \|\nabla u\|_2^2-\int_{\R^N}W(x)u^2 dx\\
=&N\int_{\R^N}\tilde{G}(u)dx\\
\geq&N\frac{\alpha-2}{2}\int_{\R^N}G(u)dx.
\end{align*}
and thus
\begin{align*}
J[u]=&\frac{1}{2}\|\nabla u\|_2^2+\frac{1}{2}\int_{\R^N}V(x)u^2dx-\int_{\R^N}G(u)dx\\
\geq&\frac{1}{2}(1-\sigma_1)\|\nabla u\|_2^2-\int_{\R^N}G(u)dx\\
\geq&\left[\frac{1}{2}(1-\sigma_1)-\frac{2(1+\sigma_2)}{N(\alpha-2)}\right]\|\nabla u\|_2^2.
\end{align*}
By $(V1)-(V2)$, one can see that $\frac{1}{2}(1-\sigma_1)-\frac{2(1+\sigma_2)}{N(\alpha-2)}>0$. Hence, by Lemma \ref{lemma:20210611-zl1}, we have that $C_a>0$.
\er

\bc\lab{cro:20210617-c1}
Under the assumptions (G1)-(G3) and $(V1)$-$(V3)$, $J\big|_{\mathcal{P}_{a}}$ is coercive , i.e.,
$$\lim_{u\in \mathcal{P}_{a}, \|\nabla u\|\rightarrow \infty}J[u]=+\infty.$$
\ec
\bp
By Remark \ref{remark:Ca-positive}, for $u\in \mathcal{P}_{a}$, we have that
$$J[u]\geq \left[\frac{1}{2}(1-\sigma_1)-\frac{2(1+\sigma_2)}{N(\alpha-2)}\right]\|\nabla u\|_2^2.$$
Hence, $$\lim_{u\in \mathcal{P}_{a}, \|\nabla u\|\rightarrow \infty}J[u]=+\infty.$$
\ep

\s{Proof of Theorem \ref{thm:non-potential}}\lab{Proof-th1}
\renewcommand{\theequation}{4.\arabic{equation}}
\bp
Let $\{u_n\}\subset \mathcal{P}_a$ such that $J[u_n]\rightarrow C_a>0$ (see Remark \ref{remark:Ca-positive}). On the other hand, we can prove that $C_a<m_a$. Indeed, let $\omega_a\in S_a$ attains $m_a$, by maximum principle, one can see that $\omega_a(x)>0$ in $\R^N$.
Since $0\not\equiv V(x)$ and $\sup_{x\in \R^N}V(x)=0$, it is easy to see that
\begin{align}
C_a\leq& \max_{t>0}J[t\star \omega_a]=J[t_{\omega_a}\star \omega_a]=I[t_{\omega_a}\star \omega_a]+\frac{1}{2}\int_{\R^N}V(x)\omega_a^2(x)dx\nonumber\\
<&I[t_{\omega_a}\star \omega_a]\leq \max_{t>0}I[t\star \omega_a]=I[\omega_a]=m_a.
\end{align}

By Corollary \ref{cro:20210617-c1}, $\{u_n\}$ is bounded in $H^1(\R^N)$. Up to a subsequence, we may assume that $u_n\rightharpoonup u$ in $H^1(\R^N)$. Then we have $u\neq 0$. If not, by the well known Brezis-Lieb lemma,
$C_a+o(1)=I[u_n]$ and $\Psi'_{\infty, u_n}(1)=o(1)$. Then by the uniqueness, there exists $t_n=1+o(1)$ such that $t_n\star u_n\in \mathcal{P}_{\infty,a}$, (see Remark \ref{remark:for-I}-(1)).
Then it is easy to see that (see also \cite[Lemma 2.7]{Yang2020})
$$m_a\leq I[t_n\star u_n]=I[u_n]+o(1)=C_a+o(1),$$
a contradiction to the fact $C_a<m_a$.

Since $\{u_n\}$ is bounded in $H^1(\R^N)$, it is easy to see that
$$\lambda_n:=-\frac{\langle J'[u_n], u_n\rangle}{a}$$
is a bounded sequence.
Furthermore,
\begin{align*}
\lambda_n \|u_n\|_2^2=& -\langle J'[u_n], u_n\rangle\\
=&\int_{\R^N}g(u_n)u_ndx -\int_{\R^N}V(x)u_n^2dx-\|\nabla u_n\|_2^2\\
\geq&\int_{\R^N}g(u_n)u_ndx-\|\nabla u_n\|_2^2\\
=&\int_{\R^N}g(u_n)u_ndx-\left\{\int_{\R^N}W(x)u_n^2dx +N\int_{\R^N}[\frac{1}{2}g(u_n)u_n-G(u_n)]dx\right\}\\
=&N\int_{\R^N}G(u_n)dx - \frac{N-2}{2}\int_{\R^N}g(u_n)u_ndx -\int_{\R^N}W(x)u_n^2dx\\
\geq&\left[N-\frac{\beta (N-2)}{2}\right]\int_{\R^N}G(u_n)dx -\sigma_2\|\nabla u_n\|_2^2.
\end{align*}
On the other hand,
\begin{align*}
(1-\sigma_2)\|\nabla u_n\|_2^2\leq &\|\nabla u_n\|_2^2-\int_{\R^N}W(x)u_n^2dx\\
=&N\int_{\R^N}[\frac{1}{2}g(u_n)u_n-G(u_n)]dx\\
\leq& \frac{N(\beta-2)}{2}\int_{\R^N}G(u_n)dx.
\end{align*}
Hence,
$$\lambda_n \|u_n\|_2^2\geq \left\{\left[N-\frac{\beta (N-2)}{2}\right]\frac{2}{N(\beta-2)}(1-\sigma_2)-\sigma_2\right\}\|\nabla u_n\|_2^2.$$
By $(V2)$, $\sigma_2<\frac{N}{\beta}-\frac{N-2}{2}$, which is equivalent to
$$\left[N-\frac{\beta (N-2)}{2}\right]\frac{2}{N(\beta-2)}(1-\sigma_2)-\sigma_2>0.$$
Hence, by Lemma \ref{lemma:20210611-zl1}, there exists some $\delta>0$ such that
$\lambda_n a > \delta$ for all $n\in \N$.
Up to a subsequence, we may assume that $\lambda_n\rightarrow \lambda>0$. One can see that $u$ solves
$$-\Delta u+V(x)u+\lambda u=g(u)\;\hbox{in}\;\R^N.$$
Put $b:=\|u\|_2^2$, then $u\in \mathcal{P}_b$ and thus $J[u]>0$.
We claim that $b=a$. If not, $c:=a-b\in (0,a)$. Let $\phi_n:=u_n-u$, by Brezis-Lieb lemma again, one can prove that $\|\phi_n\|_2^2=c+o(1), C_a+o(1)=J[u]+I[\phi_n]$ and
$\Psi'_{\infty, \phi_n}(1)=o(1)$. We claim that $$\liminf_{n\rightarrow \infty}\|\nabla \phi_n\|_2^2>0.$$
If not, by Gagliardo-Nirenberg inequality, we also have that $\liminf_{n\rightarrow \infty}\int_{\R^N}\tilde{G}(\phi_n)=0$. And then we have that $\lambda c=0$, a contradiction.
Furthermore, by the uniqueness again (see Remark \ref{remark:for-I}-(1)), there exists $t_n=1+o(1)$ such that $t_n\star \phi_n\in \mathcal{P}_{\infty, \|\phi_n\|_2^2}$. By \cite[Lemma 2.7]{Yang2020} again,
$$m_{\|\phi_n\|_2^2}\leq I[t_n\star \phi_n]=I[\phi_n]+o(1)=C_a-J[u]+o(1).$$
Noting that $m_a$ is continuous and decreasing respect to $a\in \R^+$ (see \cite[Theorem 1.1]{Yang2020}), we have that
$$m_a<m_c\leq \lim_{n\rightarrow +\infty} m_{\|\phi_n\|_2^2}\leq \lim_{n\rightarrow \infty} I[t_n\star \phi_n]=C_a-J[u]<m_a-J[u],$$
also a contradiction to $J[u]>0$.

Hence, we prove that $u_n\rightarrow u$ in $L^2(\R^N)$ and thus $u\in \mathcal{P}_a$. So $J[u]\geq C_a$.
Under the assumption $(G2)$, by Gagliardo-Nirenberg  inequality, we have
$\displaystyle \int_{\R^N}G(u_n)dx\rightarrow \int_{\R^N}G(u)dx$ as $n\rightarrow \infty$.
Furthermore,
$\displaystyle\|\nabla u\|_2^2\leq \liminf_{n\rightarrow \infty}\|\nabla u_n\|_2^2$. So we also have that
$$J[u]\leq \liminf_{n\rightarrow \infty}J[u_n]=C_a.$$
Hence, $J[u]=C_a$ and $\|\nabla u\|_2^2=\lim_{n\rightarrow \infty}\|\nabla u_n\|_2^2$. That is, $u_n\rightarrow u$ in $H^1(\R^N)$ and $u\in \mathcal{P}_a$ attains $C_a$.
By Proposition \ref{prop:20210615-p1}, there exists some $\lambda\in \R$ such that
$(\lambda_a, u_a):=(\lambda, u)$ solves \eqref{eq:problem-nonpotential}.
\ep

\end{document}